\documentclass[preprint,12pt]{elsarticle}

\usepackage{amsfonts}
\usepackage{amssymb}
\usepackage{eucal}
\usepackage{amsxtra}
\usepackage{makeidx}
\usepackage{xypic}
\usepackage{mathrsfs}
\usepackage{pdflscape}
\usepackage{tabularx}
\usepackage{xy}
\usepackage{MnSymbol}

\usepackage{amsmath}
\numberwithin{equation}{section}

\newtheorem{theorem}{Theorem}[section]
\newtheorem{proposition}[theorem]{Proposition}
\newtheorem{definition}[theorem]{Definition}
\newtheorem{lemma}[theorem]{Lemma}
\newtheorem{corollary}[theorem]{Corollary}
\newtheorem{example}[theorem]{Example}

\newproof{proof}{Proof}

\journal{Journal of Algebra}

\begin{document}

\begin{frontmatter}

\title{Weak crossed product orders arising from partitions of a finite group}

\author{Ch.Lamprakis}
\ead{chrlambr@hotmail.com}

\author{Th.Theohari-Apostolidi}
\ead{theohari@math.auth.gr}

\address{School of Mathematics, Aristotle University of Thessaloniki, Thessaloniki 54124, Greece}

\begin{abstract}
Let $R$ be a complete discrete valuation ring with quotient field $K$, $L$ a finite Galois extension of $K$ with Galois group $G$ and $S$ the integral closure of $R$ in $L$. In this article, using elements of the monoid $Sl(G)$ introduced in \cite{LaTh}, we construct certain crossed product orders in the case that the extension $S/R$ is unramified. We associate an element of $Sl(G)$ to every crossed product order and give a criterion for inflating idempotent $2$-cocycles.
\end{abstract}


\end{frontmatter}

\section{Introduction}

Let $R$ be a complete discrete valuation ring with quotient field $K$, $L$ a finite Galois extension of $K$ with Galois group $G$, $S$ the integral closure of $R$ in $L$ and $\pi S$ (resp.$\pi_0 R$) the unique maximal ideal of $S$ (resp.$R$), $L^*:=L\setminus \{0\}$, $S^*:=S\setminus \{0\}$ and $U(S)$ the group of units of $S$.

A function $f:G\times G\longrightarrow L^*$ is called normalized $2$-cocycle when the following conditions are satisfied:
	\begin{eqnarray}
		\label{eqcoccond1} && f(\sigma,\tau)f(\sigma\tau,\rho)=f^{\sigma}(\tau,\rho)f(\sigma,\tau\rho),\quad \sigma,\tau\rho\in G,\\
		\label{eqcoccond2} && f(1,\sigma)=f(\sigma,1)=1,\quad\text{for}\;\sigma\in G,
	\end{eqnarray}
where $l^{\sigma}$ means $\sigma(l)$, for $\sigma\in G$ and $l\in L$, and $1$ is the unit element of $G$. Associated to a $2$-cocycle $f$ there is a $K$-algebra $A_f$ called the crossed product algebra associated to $f$. The $K$-algebra $A_f$ is defined as an $L$-vector space
	\begin{displaymath}
		A_f=\bigoplus_{\sigma\in G} Lx_{\sigma}
	\end{displaymath}
having as an $L$-basis the symbols $x_{\sigma}$, $\sigma\in G$, and multiplication defined by the rules
	\begin{displaymath}
		x_{\sigma}l=l^{\sigma}x_{\sigma}\quad\text{and}\quad x_{\sigma}x_{\tau}=f(\sigma,\tau)x_{\sigma\tau},
	\end{displaymath}
for $\sigma,\tau\in G$ and $l\in L$. The cocycle condition (\ref{eqcoccond1}) guarantees that $A_f$ is an associative $K$-algebra with unit element $x_1$, that we denote also by $1$, and moreover $A_f$ is a central simple $K$-algebra (\cite{Re},Theorem 29.6).

The set $Z^2(G,L^*)$ of all $2$-cocycles of $G$ on $L^*$ becomes a group under the pointwise multiplication called the group of boundaries of $G$ on $L^*$ (\cite{DF}). The subgroup $B^2(G,L^*)$ of $Z^2(G,L^*)$ consisting of $2$-cocycles $b:G\times G\longrightarrow L^*$ such that
	\begin{displaymath}
		b(\sigma,\tau)=a(\sigma)a^{\sigma}(\tau)a(\sigma\tau)^{-1},
	\end{displaymath}
for $\sigma,\tau\in G$ and for a function $a:G\longrightarrow L^*$, is called the group of coboundaries of $G$ on $L^*$ and the quotient group $Z^2(G,L^*)/B^2(G,L^*)$ is the second cohomology group $H^2(G,L)$ (\cite{Re},Section 29).

If we allow $f$ to take values on $L$ instead of $L^*$, then we are reffering to a weak $2$-cocycle of $G$ on $L$. D.E. Haile, R.G. Larson and M.E. Sweedler in \cite{HLS} studied the weak $2$-cocycles and they introduced a new cohomology theory based on them. The set $Z^2(G,L)$ of $2$-cocycles of $G$ on $L$ becomes a monoid under the pointwise multiplication and
	\begin{displaymath}
		M^2(G,L):=Z^2(G,L)/ B^2(G,L^*)
	\end{displaymath}
is also a monoid, the invertible elements of which are exactly the elements of $H^2(G,L)$. Moreover every class $[f]\in M^2(G,L)$ contains a unique idempotent $2$-cocycle, i.e. that takes values on the set $\{0,1\}$ of the idempotent elements of $L$.

A function $f:G\times G\longrightarrow S^*$ satisfying relations (\ref{eqcoccond1}) and (\ref{eqcoccond2}) is called a normalized weak $2$-cocycle of $G$ on $S^*$. The ring
	\begin{displaymath}
		\Lambda_f:=\bigoplus_{\sigma\in G} Sx_{\sigma}:=(S/R,f)
	\end{displaymath}
is an $R$-order in the crossed product $K$-algebra $A_f$ and is called the weak crossed product order associated to $f$. D.E. Haile introduced and studied the weak crossed product order $\Lambda_f$ in case $R$ is a discrete valuation ring and the extension $L/K$ is unramified (\cite{Ha1}). Let $Z^2(G,S^*)$ be the set of boundaries of $G$ on $S^*$, that is, the set of all weak $2$-cocycles of $G$ on $S^*$. Then $Z^2(G,S^*)$ becomes a monoid under the pointwise multiplication having as a subgroup the set $B^2(G,U(S))$ consisting of all functions $b:G\times G\longrightarrow U(S)$ such that
	\begin{displaymath}
		b(\sigma,\tau)=a(\sigma)a^{\sigma}(\tau)a(\sigma\tau)^{-1},
	\end{displaymath}
for $\sigma,\tau\in G$ and a function $a:G\longrightarrow U(S)$, called the group of coboundaries of $G$ on $U(S)$. The set $Z^2(G,S^*)/B^2(G,U(S))$ of equivalent classes, denoted $N^2(G,S)$, is a monoid under the pointwise multiplication. Since every element of $Z^2(G,L)$ or of $Z^2(G,S^*)$ is cohomologous to  a normalized $2$-cocycle, in the whole of the article we always take normalized $2$-cocycles.

In this article we assume that the extension $L/K$ is unramified with finite residue field. Then the extension $S/R$ is also a Galois extension and $G\cong Gal(S/R)\cong Gal(\overline{S}/\overline{R})$, where $\overline{S}:=S/\pi S$ and $\overline{R}:=R/\pi_0 R$.

In \cite{LaTh}, the authors introduced the monoid $Sl(G)$, for a finite group $G$, consisting of all functions $r:G\longrightarrow \mathbb{N}$ satisfying the relations $r(1)=0$ and $r(\sigma\tau)\leq r(\sigma)+r(\tau)$, for all $\sigma,\tau\in G$. These functions arose naturally from some partitions of $G$ and were used to construct idempotent $2$-cocycles by the rule 
	\begin{equation}{\label{eqidcocyle}}
		e_r(\sigma,\tau)=
			\begin{cases}
				1, & \text{if}\; r(\sigma\tau)=r(\sigma)+r(\tau),\\
				0, & \text{if}\; r(\sigma\tau)<r(\sigma)+r(\tau).
			\end{cases}
	\end{equation}
(\cite{LaTh},Theorem 5.2).

In this article, using the elements of $Sl(G)$ we establish a relation between weak crossed product orders of $G$ on $S^*$ and weak crossed product algebras of $G$ on $L$. It consists of five sections including Introduction. In Section 2 we construct new elements of $Sl(G)$ from a given one. In Section 3 we define the monoid $B^2(G,S^*)\subseteq Z^2(G,S^*)$ and we prove that the weak crossed product orders associated to elements of this set are not hereditary. Moreover we prove that the monoid $Sl(G)$ is all we need to describe up to equivalence the monoid $B^2(G,S^*)$. In Section 4 we associate an element of $Sl(G)$ to a weak $2$-cocycle $f:G\times G\longrightarrow S^*$ and we give two applications. Finally in Section 5, for a normal subgroup $H$ of $G$, we give a necessary and sufficient condition for inflating idempotent weak $2$-cocycles from $G/H$ to $G$ and specialize this procedure using elements of $Sl(G)$.

We refer to \cite{HLS} and \cite{Re} for Galois theory and Galois cohomology theory of groups, to \cite{Re} for crossed product algebras, to \cite{CR} for orders and crossed product orders, to \cite{Ha1} and \cite{Ke} for weak crossed product orders, to \cite{Se} for local fields and ramification theory and to \cite{CR}, \cite{Lam}, \cite{Pi} for ring theory and algebras.

\section{Preliminaries}

Throughout the rest of the article we consider the notation and the assumptions of Section 1. Moreover we denote by $E^2(G,L)$ the set of idempotent $2$-cocycles of $G$ on $L$.

One of the main results of (\cite{HLS},Section 7) is that
	\begin{displaymath}
		M^2(G,L)=\bigcupdot M_e^2(G,L),
	\end{displaymath}
where the union is over all idempotent $2$-cocycles $e$ and
	\begin{eqnarray}
		M_e^2(G,L):&=& \{[f]\in M^2(G,L):[f][e]=[f]\;\text{and}\exists [g]\in M^2(G,L)\nonumber\\
			&& \text{such that}\;[f][g]=[e]\},\nonumber
	\end{eqnarray}
which is a group with unit the class $[e]$. For an element $f\in Z^2(G,L)$, let $H(f)=\{\sigma\in G:f(\sigma,\sigma^{-1})\neq 0\}$. It was shown in (\cite{HLS},Section 10) that $H(f)$ is a subgroup of $G$, called the inertial group of $f$. Moreover $H(f)=H(e)$, where $e$ is the unique idempotent $2$-cocycle such that $[f]\in M_e^2(G,L)$. Similarly if $f\in Z^2(G,S^*)$, then the set $\{\sigma\in G:f(\sigma,\sigma^{-1})\in U(S)\}$ is also a subgroup of $G$ (\cite{Ha2}). We call this subgroup also the inertial group of the $2$-cocycle $f$ and we denote it also as $H(f)$. Occasionally we use the notion $H_G(f)$ for $H(f)$ to emphasise the group $G$.

The set $Sl(G)$ becomes an additive monoid with operation $(r_1+r_2)(\sigma)=r_1(\sigma)+r_2(\sigma)$, for $r_1,r_2\in Sl(G)$ and $\sigma\in G$, with unit element the function $r(\sigma)=0$, for all $\sigma\in G$ (\cite{LaTh},Proposition 8.1). Moreover the set $M_r:=\{\sigma\in G:r(\sigma)=0\}$ is a subgroup of $G$ (\cite{LaTh},Proposition 4.2).

The next proposition is a criterion for checking the equality of elements of $Sl(G)$.

\begin{proposition}{\label{propbronto}}
Let $r_1,r_2\in Sl(G)$. Then $r_1=r_2$ if and only if
	\begin{displaymath}
		r_1(\sigma)+r_1(\tau)-r_1(\sigma \tau)=r_2(\sigma)+r_2(\tau)-r_2(\sigma \tau),
	\end{displaymath}
for $\sigma,\tau\in G$.
\end{proposition}
\textbf{Proof}: If $r_1=r_2$, then the equality is obvious. Now let $r_1(\sigma)+r_1(\tau)-r_1(\sigma \tau)=r_2(\sigma)+r_2(\tau)-r_2(\sigma \tau)$, for $\sigma,\tau\in G$. First we show that $M_{r_1}=M_{r_2}$. Both sets are not empty since they contain the unit of $G$. Let $\sigma\in M_{r_1}$. Then, for $\tau=\sigma^{-1}$, we have $r_1(\sigma)+r_1(\sigma^{-1})-r_1(1)=r_2(\sigma)+r_2(\sigma^{-1})-r_2(1)$. Since $r_1,r_2\in Sl(G)$ it follows that $r_1(1)=r_2(1)=0$. Also $\sigma^{-1}\in M_{r_1}$ so $r_1(\sigma^{-1})=0$. We deduce that $r_2(\sigma)+r_2(\sigma^{-1})=0$ and since both numbers in the left side of this equation are positive it follows that $r_2(\sigma)=r_2(\sigma^{-1})=0$, from where $\sigma\in M_{r_2}$. So $M_{r_1}\subseteq M_{r_2}$. With the same argument we prove that $M_{r_2}\subseteq M_{r_1}$ and so $M_{r_1}=M_{r_2}$.

Suppose now that $r_1\neq r_2$. We set $M=M_{r_1}=M_{r_2}$. Then there exists $\rho\in G\setminus M$ for which $r_1(\rho)\neq r_2(\rho)$ (of course $\rho\notin M)$. Suppose that $r_1(\rho)<r_2(\rho)$ (otherwise we interchange the roles of $r_1$ and $r_2$). Firstly we show by induction that $r_1(\rho^k)<r_2(\rho^k)$, for every $k\in \mathbb{N}^*:=\mathbb{N}\setminus\{0\}$. For $k=1$, we have nothing to prove. We suppose the above relation holds for some $k\in\mathbb{N}^*$. By the assumption, setting $\sigma=\rho$ and $\tau=\rho^k$, we have $r_1(\rho)+r_1(\rho^k)-r_1(\rho^{k+1})=r_2(\rho)+r_2(\rho^k)-r_2(\rho^{k+1})$. From the induction hypothesis, $r_1(\rho^k)<r_2(\rho^k)$ and also $r_1(\rho)<r_2(\rho)$, so
	\begin{displaymath}
		r_2(\rho^{k+1})-r_1(\rho^{k+1})=(r_2(\rho^k)-r_1(\rho^k))+(r_2(\rho)-r_1(\rho))>0,
	\end{displaymath}
therefore $r_1(\rho^{k+1})<r_2(\rho^{k+1})$ and the induction is completed.

The group $G$ is finite so, for $\rho\in G\setminus \{1\}$, $ord(\rho)-1$ belongs to $\mathbb{N}^*$ and then $r_1(\rho^{ord(\rho)-1})<r_2(\rho^{ord(\rho)-1})$. On the other hand by the assumption, setting $\sigma=\rho$ and $\tau=\rho^{-1}$ we have
	\begin{eqnarray}
		&& r_1(\rho)+r_1(\rho^{-1})-r_1(\rho\rho^{-1})=r_2(\rho)+r_2(\rho^{-1})-r_2(\rho\rho^{-1})\Rightarrow\nonumber\\
		&& r_1(\rho^{-1})-r_2(\rho^{-1})=r_2(\rho)-r_1(\rho)>0\Rightarrow r_1(\rho^{ord(\rho)-1})>r_2(\rho^{ord(\rho)-1}),\nonumber
	\end{eqnarray}
which is a contradiction. Hence $r_1=r_2$.\quad$\square$\\

Let $e$ be an idempotent $2$-cocycle and $H=H(e)$, then associated to $e$ is a partial ordering on the $G$-set $G/H$ given by
	\begin{displaymath}
		\sigma H\leq \tau H\Leftrightarrow e(\sigma,\sigma^{-1}\tau)=1,
	\end{displaymath}
for $\sigma,\tau\in G$ (\cite{HLS},Section 7). It is shown in (\cite{HLS},Theorem 7.13) that there is a one-to-one correspondence between the idempotent $2$-cocycles with inertial group $H$ and the partial orderings on $G/H$, with least element $H$ which are  \textquotedblleft lower subtractive\textquotedblright, that is, satisfy
	\begin{displaymath}
		\sigma H\leq \tau H\Rightarrow (\sigma H\leq \rho H\leq \tau H\Leftrightarrow \sigma^{-1}\rho H\leq \sigma^{-1}\tau H),
	\end{displaymath}
for $\sigma,\tau,\rho\in G$. The partial ordering associated to the idempotent $e$ gives rise to a left graph with vertices the left cosets of $H$ in $G$ in the obvious way (\cite{HLS},Section 7). Let $r\in Sl(G)$, $e_r$ as in relation (\ref{eqidcocyle}) and $H=H(e_r)$, then $M_r=H$ (\cite{LaTh},Theorem 5.2) and
	\begin{displaymath}
		\sigma H\leq \tau H\Leftrightarrow r(\tau)=r(\sigma)r(\sigma^{-1}\tau),
	\end{displaymath}
for $\sigma,\tau\in G$ (\cite{LaTh},Remark 5.3). This relation allows us to construct the corresponding graphs on the left cosets.

By the following proposition we can construct new elements of $Sl(G)$ from a given one.

\begin{proposition}{\label{propracc}}
Let $r\in Sl(G)$, $e_r$ be as in relation (\ref{eqidcocyle}), $H=H(e_r)$ and
	\begin{displaymath}
		N_1(e_r)=\{\sigma\in G:u_{\sigma}\in rad(A_{e_r})\;\text{and}\;u_{\sigma}\notin rad(A_{e_r})^2\}.
	\end{displaymath} 
Then, for an element $a\in N_1(e_r)$, the function $r'_a:G\rightarrow \mathbb{N}$ defined by the rule
	\begin{displaymath}
		r'_a(\sigma)=
			\begin{cases}
				r(\sigma), & \text{if}\; \sigma\notin HaH,\\
				r(\sigma)+1, & \text{if}\; \sigma\in HaH,
			\end{cases}
	\end{displaymath}
is an element of $Sl(G)$ with $M_{r'_a}=H$.
\end{proposition}
\textbf{Proof}: Let $r'_a$ be as in the statement of the proposition and let us denote $r'_a=r'$. We remark that from the definition of $r'$ it holds that $r(\sigma)\leq r'(\sigma)$, for $\sigma\in G$. From (\cite{LaTh},Lemma 3.1) we see that $a\notin H$ and consequently $1\notin HaH$. Hence $r'(1)=r(1)=0$. Now, for $\sigma,\tau\in G$, if $\sigma\tau\notin HaH$, then $r'(\sigma\tau)=r(\sigma\tau)\leq r(\sigma)+r(\tau)\leq r'(\sigma)+r'(\tau)$. So it remains to examine the case $\sigma\tau\in HaH$, for $\sigma,\tau\in G$. In this case $r'(\sigma\tau)=r(\sigma\tau)+1$. We have the following three subcases: i) $\sigma\in H$, ii) $\tau\in H$ and iii) $\sigma,\tau\notin H$. If $\sigma\in H$, then $\tau\in HaH$ and so $r'(\sigma)=r(\sigma)$ and $r'(\tau)=r(\tau)+1$. Hence $r'(\sigma\tau)\leq r'(\sigma)+r'(\tau)$. If $\tau\in H$, then similarly we get $r'(\sigma\tau)\leq r'(\sigma)+r'(\tau)$. Finally, if $\sigma,\tau\notin H$, then, since $HaH\subseteq  N_1(e_r)$ (\cite{LaTh},Remark 3.2) and from (\cite{LaTh},Lemma 3.1) we get $e_r(\sigma,\tau)=0$ and hence $r(\sigma\tau)<r(\sigma)+r(\tau)$. Therefore $r'(\sigma\tau)=r(\sigma\tau)+1<r(\sigma)+r(\tau)+1\leq r'(\sigma)+r'(\tau)+1$ and since all terms are natural numbers we get $r'(\sigma\tau)\leq r'(\sigma)+r'(\tau)$. So we conclude that $r'\in Sl(G)$.

To prove that $H\subseteq M_{r'}$ we remark that if $\sigma\in H$, then $\sigma\notin HaH$ and hence $r'(\sigma)=r(\sigma)=0$. For the other direction if $\sigma\in M_{r'}$, then $r'(\sigma)=0$ and since $r(\sigma)\leq r'(\sigma)$, for all $\sigma\in G$, we get that $r(\sigma)=0$ and so $M_{r'}\subseteq H$.\quad$\square$

\section{Constructing weak crossed product orders from elements of $Sl(G)$}

We consider the $\pi S$-adic valuation $v_{\pi}:S\longrightarrow \mathbb{Z}\cup \{\infty\}$, for $v_{\pi}(s)=max\{z\in\mathbb{Z}:s\in\pi^z S\}$. Let $\sigma(\pi)=u(\sigma)\pi$, where $u(\sigma)\in U(S)$, for $\sigma\in G$, and we define the function $u:G\longrightarrow U(S)$, $\sigma\mapsto u(\sigma)$. We refer to (\cite{Se}) for more details.

\begin{definition}{\label{defcobsstar}}
We define $B^2(G,S^*)$ to be the set of all functions $b:G\times G\longrightarrow S^*$ such that there exists a function $a:G\longrightarrow S^*$ with $a(1)=1$, for which
	\begin{displaymath}
		b(\sigma,\tau)=a(\sigma)a^{\sigma}(\tau)a(\sigma \tau)^{-1},\quad \text{for}\;\sigma,\tau\in G.
	\end{displaymath}
\end{definition}
>From the definition of $B^2(G,S^*)$ it is clear that $a(\sigma\tau)^{-1}\in L^*$ and so $b(\sigma,\tau)$ is not in general an element of $S^*$. However, we are interested only in functions $a$ forcing $b$ to take values in $S^*$. For example, if the function $a$ takes values in $U(S)\subseteq S^*$, then the function $b$ takes values in $U(S)$ and then $b$ is an element of $B^2(G,U(S))$. This demonstrates that
	\begin{displaymath}
		B^2(G,U(S))\subseteq B^2(G,S^*)\subseteq B^2(G,L^*).
	\end{displaymath}
It is easy to show that the set $B^2(G,S^*)$ is a submonoid of $Z^2(G,S^*)$.

\begin{proposition}{\label{propeqbr}}
Let $r\in Sl(G)$, then the function $b_r:G\times G\longrightarrow S^*$ defined by then rule 
	\begin{equation}{\label{eqbr}}
		b_r(\sigma,\tau)=u(\sigma)^{r(\tau)}\pi^{r(\sigma)+r(\tau)-r(\sigma \tau)},\quad\text{for}\;\sigma,\tau\in G,
	\end{equation}	
is an element of $B^2(G,S^*)$ with inertial group $M_r$.
\end{proposition}
\textbf{Proof}: For any $r\in Sl(G)$ consider the function $a_r:G\longrightarrow S^*$ defined by the rule $a_r(\sigma)=\pi^{r(\sigma)}$. Then $a_r(\sigma)a_r^{\sigma}(\tau)a_r(\sigma \tau)^{-1}=b_r(\sigma,\tau)$. Since $r\in Sl(G)$ and so $r(\sigma)+r(\tau)\geq r(\sigma \tau)$, it follows that $b_r(\sigma,\tau)$ is a $2$-cocycle taking values in $S^*$, for $\sigma,\tau\in G$, so $b_r\in B^2(G,S^*)$. Finally for the inertial group $H(b_r)$ of $b_r$ we have
	\begin{displaymath}
		H(b_r)= \{\sigma\in G:b_r(\sigma,\sigma^{-1})\in U(S)\}=\{\sigma\in G:r(\sigma)+r(\sigma^{-1})=0\}=M_r.\nonumber\square
	\end{displaymath}

In what follows, for $r\in Sl(G)$, we use the notation
	\begin{equation}{\label{eqepsilonr}}
		\epsilon_r(\sigma,\tau)=\pi^{r(\sigma)+r(\tau)-r(\sigma \tau)},
	\end{equation}
so that relation (\ref{eqbr}) takes the form
	\begin{equation}{\label{eqbrepsilonr}}
		b_r(\sigma,\tau)=u(\sigma)^{r(\tau)}\epsilon_r(\sigma,\tau).
	\end{equation}
The weak $2$-cocycle $b_r$ defines a weak crossed product order $\Lambda_{b_r}=(S/R,b_r)$ and the idempotent $2$-cocycle $e_r$ of relation (\ref{eqidcocyle}) defines a weak crossed product algebra $A_{e_r}=(L/K,e_r)$. By reducing $mod\pi S$ the values of $b_r$ of relation (\ref{eqbr}) we get a new weak $2$-cocycle $\overline{b}_r$ of $G$ in $\overline{S}$. More explicitly,
	\begin{displaymath}
		\overline{b}_r(\sigma,\tau)=
			\begin{cases}
				u(\sigma)^{r(\tau)}+\pi S, & \text{if}\; r(\sigma\tau)=r(\sigma)+r(\tau),\\
				0, & \text{if}\; r(\sigma\tau)<r(\sigma)+r(\tau).
			\end{cases}
	\end{displaymath}

\begin{example}{\label{exbr}}
\end{example}
Let $G=\mathbb{Z}/10\mathbb{Z}$ be the cyclic group of order $10$. Consider the standard partition of $G\setminus \{\overline{0}\}$ derived from the generating set $\{\overline{1},\overline{6}\}$ (\cite{Al1},Section 3 and \cite{LaTh},Section 10):
	\begin{displaymath}
		\delta= \{\{\overline{1},\overline{6}\},\{\overline{2},\overline{7}\},\{\overline{3},\overline{8}\},\{\overline{4},\overline{9}\},\{\overline{5}\}\},
	\end{displaymath}
from where we can construct $r=\{0,1,2,3,4,5,1,2,3,4\}\in Sl(G)$ with $N_1(e_r)=\{\overline{1},\overline{6}\}$ (\cite{LaTh},Remark 10.9 ). With the procedure described in Proposition \ref{propracc}, we can construct the function
	\begin{displaymath}
		r_1=r'_{\overline{1}}=\{0,2,2,3,4,5,1,2,3,4\},
	\end{displaymath}
which is again an element of $Sl(G)$. From relation (\ref{eqepsilonr}), the values of the non-invertible term $\epsilon_{r_1}:G\times G\longrightarrow S^*$ are given by the following table:
	\begin{displaymath}
		\begin{array}{|c|cccccccccc|}
		\hline
		(\sigma,\tau)&\overline{0}&\overline{1}&\overline{2}&\overline{3}&\overline{4}&
\overline{5}&\overline{6}&\overline{7}&\overline{8}&\overline{9}\\
		\hline
			\overline{0}&1&1&1&1&1&1&1&1&1&1\\
			\overline{1}&1&\pi^2&\pi&\pi&\pi&\pi^6&\pi&\pi&\pi&\pi^6\\
			\overline{2}&1&\pi&1&1&\pi^5&\pi^5&1&1&\pi^5&\pi^4\\
			\overline{3}&1&\pi&1&\pi^5&\pi^5&\pi^5&1&\pi^5&\pi^4&\pi^5\\
			\overline{4}&1&\pi&\pi^5&\pi^5&\pi^5&\pi^5&\pi^5&\pi^4&\pi^5&\pi^5\\
			\overline{5}&1&\pi^6&\pi^5&\pi^5&\pi^5&\pi^{10}&\pi^4&\pi^5&\pi^5&\pi^5\\
			\overline{6}&1&\pi&1&1&\pi^5&\pi^4&1&1&1&1\\
			\overline{7}&1&\pi&1&\pi^5&\pi^4&\pi^5&1&1&1&\pi^5\\
			\overline{8}&1&\pi&\pi^5&\pi^4&\pi^5&\pi^5&1&1&\pi^5&\pi^5\\
			\overline{9}&1&\pi^6&\pi^4&\pi^5&\pi^5&\pi^5&1&\pi^5&\pi^5&\pi^5\\
			\hline
		\end{array}
	\end{displaymath}
As in the relation (\ref{eqbrepsilonr}), multiplying the values of $\epsilon_{r_1}(\sigma,\tau)$ with the invertible term $u(\sigma)^{r_1(\tau)}$, we get the $2$-cocycle $b_{r_1}$ which has trivial inertial group.

Notice that $r_1$ defines also an idempotent $2$-cocycle $e_{r_1}$ as in relation (\ref{eqidcocyle}) with the following graph on (left) cosets:

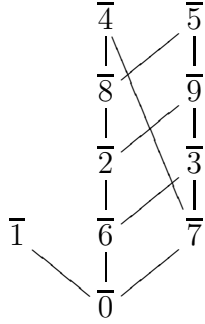
\[
\xymatrixrowsep{0.15in}
\xymatrixcolsep{0.3in}
\xymatrix{ & 	\overline{4}	&\overline{5} \\
		& \overline{8}\ar@{-}[u]\ar@{-}[ur]&\overline{9}\ar@{-}[u]\\
				& \overline{2}\ar@{-}[u]\ar@{-}[ur]&\overline{3}\ar@{-}[u]\\
				\overline{1}&\overline{6}\ar@{-}[u]\ar@{-}[ur] & \overline{7}\ar@{-}[u]\ar@{-}[uuul]\\
					& \overline{0}\ar@{-}[ul]\ar@{-}[u]\ar@{-}[ur] &}
\]

\begin{proposition}{\label{propbbr}}
Every element of $B^2(G,S^*)$ with inertial group $H$ can be written uniquely as a product of an element of $B^2(G,U(S))$ and a $2$-cocycle $b_r$, for some $r\in Sl(G)$ with $M_r=H$.
\end{proposition}
\textbf{Proof}: Let $b\in B^2(G,S^*)$. By the definition of $B^2(G,S^*)$ there exists $a:G\longrightarrow S^*$ defined by $a(\sigma)=w(\sigma)\pi^{v_{\pi}(a(\sigma))}$, with $w(\sigma)\in U(S)$ for $\sigma\in G$ and $a(1)=1$, such that $b(\sigma,\tau)=a(\sigma)\sigma(a(\tau))a(\sigma \tau)^{-1}$, for $\sigma,\tau\in G$. Substituting $a$ we get:
	\begin{displaymath}
		b(\sigma,\tau)=w(\sigma)\sigma(w(\tau))w(\sigma \tau)^{-1}u(\sigma)^{v_{\pi}(a(\tau))}\pi^{v_{\pi}(a(\sigma))+v_{\pi}(a(\tau))-v_{\pi}(a(\sigma \tau))}.
	\end{displaymath}
We set $c(\sigma,\tau)=w(\sigma)\sigma(w(\tau))w(\sigma \tau)^{-1}$ and $r=v_{\pi}\circ a$. It is obvious that $c\in B^2(G,U(S))$. Also since $b(\sigma,\tau)\in S^*$, for $\sigma,\tau\in G$, it must be $r(\sigma)+r(\tau)-r(\sigma \tau)\geq 0$. Since $r(1)=v_{\pi}(a(1))=0$ we have $r\in Sl(G)$. So $b(\sigma,\tau)=c(\sigma,\tau)b_r(\sigma,\tau)$, with $c\in B^2(G,U(S))$ and $r\in Sl(G)$.

Now, let $b(\sigma,\tau)=c'(\sigma,\tau)b_{r'}(\sigma,\tau)$ be another expression of $b$ with $c'\in B^2(G,U(S))$ and $r'\in Sl(G)$. Since the expression of $b(\sigma,\tau)$ as a product of a unit and a power of $\pi$ is unique, we must have that $r(\sigma)+r(\tau)-r(\sigma\tau)=r'(\sigma)+r'(\tau)-r'(\sigma\tau)$, for every $\sigma,\tau\in G$, and $c(\sigma,\tau)u(\sigma)^{r(\tau)}=c'(\sigma,\tau)u(\sigma)^{r'(\tau)}$. From Proposition \ref{propbronto} we deduce that $r=r'$ and as a consequence we get $c=c'$.

Finally, let $\sigma\in M_r$. Then $r(\sigma)=r(\sigma^{-1})=0$ and so $b(\sigma,\sigma^{-1})=c(\sigma,\sigma^{-1})u(\sigma)^{r(\sigma^{-1})}\pi^0\in U(S)$ from where we conclude that $\sigma\in H(b)=H$. On the other hand, if $\sigma\in H$, then $b(\sigma,\sigma^{-1})\in U(S)$ from where we deduce that $r(\sigma)+r(\sigma^{-1})=0$ and so $r(\sigma)=0$. Hence $H=M_r$.\quad$\square$

\begin{corollary}{\label{corbhered}}
Let $b\in B^2(G,S^*)$ and $H(b)\subsetneqq G$ be the inertial group of $b$. Then the weak crossed product order $\Lambda_b$ is not hereditary.
\end{corollary}
\textbf{Proof}: Let $b\in B^2(G,S^*)$ and let $b=cb_r$ be the unique expression of $b$ in accordance to Proposition \ref{propbbr}, for some $c\in B^2(G,U(S))$ and a weak $2$-cocycle $b_r:G\times G\longrightarrow S^*$. Then we have $v_{\pi}(b(\sigma,\sigma^{-1})=r(\sigma)+r(\sigma^{-1})$, for $\sigma\in G$. Since $H(b)\subsetneqq G$, let $\tau\in G\setminus H$. Then $r(\tau)\geq 1$ and $r(\tau^{-1})\geq 1$, so $r(\tau)+r(\tau^{-1})\geq 2$. Now the result follows from (\cite{Wi1},Theorem 2.15(1)) by which $\Lambda_b$ is hereditary if and only if $v_{\pi}(b(\sigma,\sigma^{-1}))\leq 1$ for every $\sigma\in G$.\quad$\square$\\

Inside the monoid $B^2(G,S^*)$ we consider the equivalence relation
	\begin{displaymath}
		b_1\sim b_2 \Leftrightarrow b_1B^2(G,U(S))=b_2B^2(G,U(S)).
	\end{displaymath}
Let $\overline{b}$ denote the class of $b$. We write $B^2(G,S^*)/B^2(G,U(S))$ for the set of classes of elements of $B^2(G,S^*)$ which is a monoid under the operation $\overline{b_1}\cdot \overline{b_2}=\overline{b_1b_2}$.

\begin{theorem}{\label{theB2eqR}}
There is a monoid isomorphism
	\begin{displaymath}
		B^2(G,S^*)/B^2(G,U(S))\simeq Sl(G)
	\end{displaymath}
\end{theorem}
\textbf{Proof}: Define the mapping $\phi$ by the rule $\phi(\overline{b})=r$, for all $b\in B^2(G,S^*)$ where $b=c\cdot b_r$ with $c\in B^2(G,U(S))$ and $r\in Sl(G)$ as in Proposition \ref{propbbr}. The uniqueness of Proposition \ref{propbbr} implies that $\phi$ is well defined. Also, if $\phi(\overline{b_1})=\phi(\overline{b_2})$, then $r_1=r_2$ and so $b_1=c_1b_r$, $b_2=c_2b_r$ from where it follows that $b_1=c_1c_2^{-1}c_2b_r=c_1c_2^{-1}b_2$ and $\phi$ is injective. To see that $\phi$ is a homomorphism it is enough to notice that $b_{r_1+r_2}=b_{r_1}b_{r_2}$. If we set $c_1(\sigma,\tau)=1$, for every $\sigma,\tau\in G$, and $r^*(\sigma)=0$, for every $\sigma\in G$, then we get the identity $2$-cocycle. Finally, for $r\in Sl(G)$, there exists $b=c_1b_r$ such that $\phi(\overline{b})=r$ and $\phi$ is onto.\quad$\square$

\section{Elements of $Sl(G)$ defined by a valuation}

In this section we construct new elements of $Sl(G)$ from elements of $Z^2(G,S^*)$.

\begin{lemma}{\label{lemvaldvr}}
Let $f\in Z^2(G,L^*)$. Then
	\begin{displaymath}
		v_{\pi}(f(\sigma,\tau))+v_{\pi}(f(\tau^{-1},\sigma^{-1}))=v_{\pi}(f(\sigma,\sigma^{-1}))+v_{\pi}(f(\tau,\tau^{-1}))-v_{\pi}(f(\sigma \tau,\tau^{-1}\sigma^{-1})).
	\end{displaymath}
\end{lemma}
\textbf{Proof}: For an element $f\in Z^2(G,L^*)$ we consider successively the cocycle condition (\ref{eqcoccond1}) for the elements $\sigma,\tau,(\sigma\tau)^{-1}$ and $\tau,\tau^{-1},\sigma^{-1}$ of $G$. We get two equations and then we substitute the value $f(\tau,\tau^{-1}\sigma^{-1})$ of the second equation in the first equation. In both sides of the resulting equation we apply the valuation $v_{\pi}$ and we get the result, using the fact that $v_{\pi}(l^{\sigma})=v_{\pi}(l)$, for $\sigma\in G$ and $l\in L$.\quad$\square$

\begin{proposition}{\label{proprvalpi}}
Let $f\in Z^2(G,S^*)$. Then the function
	\begin{displaymath}
		r_f:G\longrightarrow \mathbb{N},\quad r_f(\sigma)=v_{\pi}(f(\sigma,\sigma^{-1})),
	\end{displaymath}
is an element of $Sl(G)$ with $M_{r_f}=H(f)$.
\end{proposition}
\textbf{Proof}: First we see that $r_f(1)=0$. Now we remark that since $f(\sigma,\tau)\in S^*$, for $\sigma,\tau\in G$, then $r_f(\sigma)\geq 0$, $\sigma\in G$. So  from Lemma \ref{lemvaldvr} we get that
	\begin{displaymath}
		r_f(\sigma)+r_f(\tau)-r_f(\sigma\tau)\geq 0,
	\end{displaymath}
hence $r_f\in Sl(G)$. Finally
	\begin{displaymath}
		M_{r_f}=\{\sigma\in G:r_f(\sigma)=0\}=\{\sigma\in G:f(\sigma,\sigma^{-1})\in U(S)\}=H(f).\quad\square
	\end{displaymath}

Let $f\in Z^2(G,S^*)$, $r\in Sl(G)$ and 
	\begin{equation}{\label{eqcocinvval}}
		f(\sigma,\tau)=c(\sigma,\tau)\pi^{v_{\pi}(f(\sigma,\tau))},
	\end{equation}
for some $c(\sigma,\tau)\in U(S)$ and $\sigma,\tau\in G$. We define the function $h:G\times G\rightarrow S^*$ by the rule
	\begin{equation}{\label{eqhfunct}}
		h(\sigma,\tau)=[c(\sigma,\tau)]^{-1}u(\sigma)^{r(\tau)}\pi^{r(\sigma)+r(\tau)-r(\sigma \tau)-v_{\pi}(f(\sigma,\tau))},
	\end{equation}
for $\sigma,\tau\in G$ and $c(\sigma,\tau)$ is as in the relation (\ref{eqcocinvval}). The function $h\in Z^2(G,L^*)$. Indeed, it is easy to see that $fh=b_r$, where $b_r$ is as in Proposition \ref{propeqbr}. Hence $h=f^{-1}b_r$, where $f^{-1}(\sigma,\tau)=f(\sigma,\tau)^{-1}$. In the sequel we prove that, for $f\in Z^2(G,S^*)$, we can find an element $h\in Z^2(G,S^*)$ such that $fh\in B^2(G,S^*)$.

\begin{proposition}
Let $f\in Z^2(G,S^*)$ and $r_f$ as in Proposition \ref{proprvalpi}. Then the function $h:G\times G\longrightarrow S^*$ defined by the rule
	\begin{displaymath}
		h(\sigma,\tau)=[c(\sigma,\tau)c(\tau^{-1},\sigma^{-1})]^{-1}u(\sigma)^{r_f(\tau)}f(\tau^{-1},\sigma^{-1}),
	\end{displaymath}
for $c(\sigma,\tau)$ as in (\ref{eqcocinvval}) for $\sigma,\tau\in G$, is an element of $Z^2(G,S^*)$ such that $fh=b_{r_f}$.
\end{proposition}
\textbf{Proof}: We put in the relation (\ref{eqhfunct}) the function $r_f$ instead of $r$. Using Lemma \ref{lemvaldvr} and the relation $f(\tau^{-1},\sigma^{-1})=c(\tau^{-1},\sigma^{-1})\pi^{v_{\pi}(f(\tau^{-1},\sigma^{-1}))}$ we get
	\begin{displaymath}
		h(\sigma,\tau)=[c(\sigma,\tau)c(\tau^{-1},\sigma^{-1})]^{-1}u(\sigma)^{r_f(\tau)}f(\tau^{-1},\sigma^{-1}).
	\end{displaymath}
The result follows from the fact that $v_{\pi}(f(\tau^{-1},\sigma^{-1})\geq 0$.\quad$\square$\\

For the next application of the function $r_f$ we need a new transformation of elements of $Sl(G)$.

\begin{proposition}{\label{proprdiv}}
Let $r\in Sl(G)$. Then the function $r_/:G\longrightarrow \mathbb{N}$ defined by the rule
	\begin{displaymath}
		r_/(\sigma)=
			\begin{cases}
				k, & r(\sigma)=2k,\\
				k+1, & r(\sigma)=2k+1
			\end{cases}
	\end{displaymath}
is an element of $Sl(G)$ with $M_{r_/}=M_r$.
\end{proposition}
\textbf{Proof}: It is easy to see that $r_/(1)=0$. For $\sigma,\tau\in G$, we distinguish six cases in accordance to $r(\sigma)$, $r(\tau)$ or $r(\sigma\tau)$ is even or odd natural number. Let $r(\sigma)=2k_1$, $r(\tau)=2k_2$, $r(\sigma\tau)=2k_3$. Then $r_/(\sigma)=k_1$, $r_/(\tau)=k_2$, $r_/(\sigma\tau)=k_3$. Since $r\in Sl(G)$, we get $2k_3\leq 2k_1+2k_2$ from where $k_3\leq k_1+k_2$ and so $r_/(\sigma\tau)\leq r_/(\sigma)+r_/(\tau)$. Similarly we examine the other cases. Finally, $\sigma\in M_{r_/}\Leftrightarrow r_/(\sigma)=0\Leftrightarrow r(\sigma)=0$ so $M_{r_/}=M_r$.\quad$\square$

\begin{corollary}{\label{correven}}
Let $r\in Sl(G)$. Then the function $r':G\longrightarrow\mathbb{N}$ defined by the rule
	\begin{displaymath}
		r'(\sigma)=
			\begin{cases}
				r(\sigma), & r(\sigma)=2k,\\
				r(\sigma)+1, & r(\sigma)=2k+1
			\end{cases}
	\end{displaymath}
is an element of $Sl(G)$ with $M_{r'}=M_r$.
\end{corollary}
\textbf{Proof}: Let $r\in Sl(G)$ and $r_/\in Sl(G)$ be as in Proposition \ref{proprdiv}. We remark that if $r(\sigma)=2k$, then $r_/(\sigma)=k$ and so $r'(\sigma)=2k=2r_/(\sigma)$, for $\sigma\in G$. Similarly for $r(\sigma)=2k+1$, $\sigma\in G$. The result follows from (\cite{LaTh},Proposition 8.3).\quad$\square$\\

\begin{example}
\end{example}
We consider Example \ref{exbr} and let $r=\{0,2,2,3,4,5,1,2,3,4\}\in Sl(G)$. Then $r_/=\{0,1,1,2,2,3,1,1,2,2\}$ and $r'=\{0,2,2,4,4,6,2,2,4,4\}$. The graph on left cosets of the corresponding idempotent $2$-cocycles $e_{r_/}$ and $e_{r'}$ is common and it is the following:
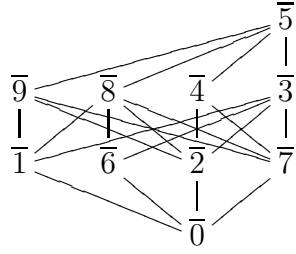
\[
\xymatrixrowsep{0.15in}
\xymatrixcolsep{0.3in}
\xymatrix{	& & &\overline{5}\\
			\overline{9}\ar@{-}[urrr]	& \overline{8}\ar@{-}[urr] & \overline{4}\ar@{-}[ur] & \overline{3}\ar@{-}[u]\\
				\overline{1}\ar@{-}[u]\ar@{-}[ur]\ar@{-}[urrr] & \overline{6}\ar@{-}[u]\ar@{-}[urr] &\overline{2}\ar@{-}[ull]\ar@{-}[ul]\ar@{-}[u]\ar@{-}[ur] & \overline{7}\ar@{-}[ulll]\ar@{-}[ull]\ar@{-}[ul]\ar@{-}[u]\\
					& &  \overline{0}\ar@{-}[ull]\ar@{-}[ul]\ar@{-}[u]\ar@{-}[ur] &}
\]

\begin{theorem}{\label{theineqordh}}
For $f\in Z^2(G,S^*)$, there exist a $2$-cocycle $h\in Z^2(G,L^*)$ and $r\in Sl(G)$ such that $fh=b_r$ with $0\leq v_{\pi}(h(\sigma,\sigma^{-1}))\leq 1$, for $\sigma\in G$.
\end{theorem}
\textbf{Proof}: Let $f\in Z^2(G,S^*)$ and $r=r_f$ as in Proposition \ref{proprvalpi}. For the $2$-cocycle $h$ of relation (\ref{eqhfunct}) and $r_/\in Sl(G)$ as in Proposition \ref{proprdiv}, we have that $fh=b_{r_/}$ and so $v_{\pi}(h(\sigma,\sigma^{-1}))=r_/(\sigma)+r_/(\sigma^{-1})-r(\sigma)$. From (\cite{Wi1},Lemma 2.6) we have that $v_{\pi}(f(\sigma,\sigma^{-1}))=v_{\pi}(f(\sigma^{-1},\sigma))$, so $r(\sigma)=r(\sigma^{-1})$ and consequently $r_/(\sigma)=r_/(\sigma^{-1})$, for $\sigma\in G$. But then
	\begin{eqnarray}
		r_/(\sigma)+r_/(\sigma^{-1})-r(\sigma) &=&
			\begin{cases}
				k+k-2k, & if\;r(\sigma)=2k,\\
				k+1+k+1-2k-1, & if\;r(\sigma)=2k+1,
			\end{cases}\nonumber\\
			&=&
			\begin{cases}
				0, & if\;r(\sigma)=2k,\\
				1, & if\;r(\sigma)=2k+1,
			\end{cases}\nonumber
	\end{eqnarray}
from where $0\leq v_{\pi}(h(\sigma,\sigma^{-1}))\leq 1$, for all $\sigma\in G$.\quad$\square$

\section{Inflating idempotent $2$-cocycles}

Let $N$ be a normal subgroup of G, $L^N$ be the fixed field of $N$ and $g\in Z^2(G/N,L^N)$. Then we can form a new normalized weak $2$-cocycle $\widehat{g}\in Z^2(G,L)$ defined by the relation
	\begin{displaymath}
		\widehat{g}(\sigma,\tau)=g(\sigma N,\tau N),\quad\text{for}\;\sigma,\tau\in G,
	\end{displaymath}
called the inflation of $g$. The main result in this section is Theorem \ref{theinertquat} by which a relation is established between the inertial groups of $\varepsilon$ and $\widehat{\varepsilon}$. First we need some lemmata. Before these let us remark that as a consequence of the definition of the inertial group $H(f)$, for $f\in Z^2(G,L)$, we get that $f(\sigma,h)\neq 0$ and $f(h,\sigma)\neq 0$, for $\sigma\in G$ and $h\in H(f)$.

\begin{lemma}{\label{lemmaus4}}
Let $f\in Z^2(G,L)$ and $N$ be a normal subgroup of $G$ such that $N\subseteq H(f)$. Then, for $\sigma,\tau\in G$, and $n_1,n_2\in N$, we have:
	\begin{itemize}
		\item[(i)] If $f(\sigma,\tau)=0$, then $f(\sigma n_1,\tau n_2)=0$,
		\item[(ii)] If $f(\sigma,\tau)\neq 0$, then \mbox{$f(\sigma n_1,\tau n_2)\neq 0$}.
	\end{itemize}
\end{lemma}
\textbf{Proof}: (i) Let $n_1,n_2\in N$. First we prove that $f(\sigma n_1,\tau)=0$. Since $N$ is a normal subgroup of $G$, there exists $n\in N$ such that $\sigma n_1=n\sigma$. From the relation (\ref{eqcoccond1}) and the elements $n,\sigma,\tau$, we have the relation $f(n,\sigma)f(n\sigma,\tau)=f^n(\sigma,\tau)f(n,\sigma \tau)$. Since $n\in H(f)$ we have that $f(n,\sigma)\neq 0$. From the assumption $f(\sigma,\tau)=0$ so $f(n\sigma,\tau)=f(\sigma n_1,\tau)=0$. Again from the relation (\ref{eqcoccond1}) and for the elements $\sigma n_1,\tau,n_2$ we have $f(\sigma n_1,\tau)f(\sigma n_1\tau,n_2)=f^{\sigma n_1}(\tau,n_2)f(\sigma n_1,\tau n_2)$. Since $f(\sigma n_1,\tau)=0$ and $f(\tau,n_2)\neq 0$, it follows that $f(\sigma n_1,\tau n_2)=0$.\\
(ii) Let $n_1,n_2\in N$. For $n\in N$ and the elements $\sigma,\tau,n$, the relation (\ref{eqcoccond1}) becomes $f(\sigma,\tau)f(\sigma \tau,n)=f^{\sigma}(\tau,n)f(\sigma,\tau n)$. Since $f(\sigma,\tau)\neq 0$ and $f(\sigma \tau,n)\neq 0$ we get $f(\sigma,\tau n)\neq 0$. Since $N$ is a normal subgroup of $G$, there exists $n'\in N$ such that $n_1\tau n_2=\tau n'$. Then from the relation (\ref{eqcoccond1}) we get
	\begin{displaymath}
		f(\sigma,n_1)f(\sigma n_1,\tau n_2)=f^{\sigma}(n_1,\tau n_2)f(\sigma,n_1\tau n_2)=f^{\sigma}(n_1,\tau n_2)f(\sigma,\tau n'),
	\end{displaymath}
since $f(n_1,\tau n_2)\neq 0$. Then $f(\sigma,\tau n')\neq 0$, so we get $f(\sigma n_1,\tau n_2)\neq 0$. $\square$\\

In what follows, for $f\in Z^2(G,L)$ and $M$ a subgroup of $G$ we write $res_M(f):M\times M\longrightarrow L$ for the restriction of $f$ to $M$, which is again a normalized weak $2$-cocycle. Also let $E^2(G,L;H)$ be the set of idempotent $2$-cocycles of $G$ in $L$ with inertial group $H$.

\begin{lemma}{\label{lemreunit}}
Let $e\in E^2(G,L;H)$ and $M$ be a subgroup of $G$. Then $res_M(e)=1$ if and only if $M\leq H$.
\end{lemma}
\textbf{Proof}: Let $e\in E^2(G,L;H)$, $e_1=res_M(e)=1$ and $\sigma\in M$. Then $\sigma^{-1}\in M$ and $e(\sigma,\sigma^{-1})=e_1(\sigma,\sigma^{-1})=1$, so $\sigma\in H$ and $M\subseteq H$. On the other hand, if $M\subseteq H$ and $\sigma\in M$, then $e_1(\sigma,\sigma^{-1})=e(\sigma,\sigma^{-1})\neq 0$ and since $e$ takes only the values 0 or 1, we get $e_1(\sigma,\sigma^{-1})=1$. Then, for $\sigma,\tau,\tau^{-1}\in M$, the $2$-cocycle condition (\ref{eqcoccond1}) yields $e_1(\sigma,\tau)e_1(\sigma \tau,\tau^{-1})=e_1(\tau,\tau^{-1})e_1(\sigma,\tau \tau^{-1})$, so $e_1(\tau,\tau^{-1})=1$. By definition $e_1(\sigma,1)=1$. Therefore $e_1(\sigma,\tau)=1$, for every $\sigma,\tau\in M$. Hence $e_1=1$.\quad$\square$\\

The following is a more general statement of (\cite{ThTo},Theorem 3.5):

\begin{theorem}{\label{theinertquat}}
Let $N$ be a normal subgroup of $G$ and $e\in E^2(G,L;H)$. Then $N\leq H$ if and only if there exists $\varepsilon\in E^2(G/N,L^N)$ such that $e=\widehat{\varepsilon}$, that is $e(\sigma,\tau)=\varepsilon(\sigma N,\tau N)$, for $\sigma,\tau\in G$. In this case
	\begin{displaymath}
		H_{G/N}(\varepsilon)=H_G(\widehat{\varepsilon})/N=H/N.
	\end{displaymath}
\end{theorem}
\textbf{Proof}: If $\sigma\in G$, then we write $\overline{\sigma}$ for $\sigma N\in G/N$. Firstly, let $N\leq H$ and $\varepsilon:G/N \times G/N\longrightarrow L^N$ defined by the rule $\varepsilon(\overline{\sigma},\overline{\tau})=e(\sigma,\tau)$, for $\overline{\sigma},\overline{\tau}\in G/N$. We prove that the function $\varepsilon$ is well defined. For this, let $\sigma,\tau\in G$, $\sigma_1=\sigma n_1$ and $\tau_1=\tau n_2$ representatives of $\overline{\sigma}$ and $\overline{\tau}$ respectively, for $n_1,n_2\in N$. Then we get the following:

i) If $e(\sigma,\tau)=1$, from Lemma \ref{lemmaus4}(ii) it follows that $e(\sigma n_1,\tau n_2)\neq 0$, so $e(\sigma_1,\tau_1)=e(\sigma n_1,\tau n_2)=1=e(\sigma,\tau )$.

ii) If $e(\sigma,\tau)=0$, from Lemma \ref{lemmaus4}(i) it follows that $e(\sigma n_1,\tau n_2)=0$, so $e(\sigma_1,\tau_1)=e(\sigma n_1,\tau n_2)=0=e(\sigma,\tau)$.

\noindent In any case $\varepsilon(\overline{\sigma}_1,\overline{\tau}_1)=e(\sigma_1,\tau_1)=e(\sigma,\tau)=\varepsilon(\overline{\sigma},\overline{\tau})$, so the value of $\varepsilon$ is independent of the representatives of $\overline{\sigma}$ and $\overline{\tau}$.

For $\overline{\sigma},\overline{\tau},\overline{\rho}\in G/N$, we get
	\begin{displaymath}
		\varepsilon(\overline{\sigma},\overline{\tau})\varepsilon(\overline{\sigma}\cdot\overline{\tau},\overline{\rho})=e(\sigma,\tau)e(\sigma \tau,\rho)=e(\tau,\rho)e(\sigma,\tau \rho)=\varepsilon(\overline{\tau},\overline{\rho})\varepsilon(\overline{\sigma},\overline{\tau}\cdot\overline{\rho}),
	\end{displaymath}
so $\varepsilon$ is an idempotent $2$-cocycle.  The $2$-cocycle $e$ is normalized, so for $\overline{\sigma},\overline{\tau}\in G/N$ it holds that $\varepsilon(\overline{1},\overline{\tau})=e(1,\tau)=1$ and $\varepsilon(\overline{\sigma},\overline{1})=e(\sigma,1)=1$. So $\varepsilon$ is also normalized.

For the opposite direction suppose that there exists $\varepsilon\in E^2(G/N,L^N)$ such that $e=\widehat{\varepsilon}$. Since $H_{G/N}(\varepsilon)\leq G/N$, there exists a subgroup $G_1$ of $G$ containing $N$ such that $H_{G/N}(\varepsilon)=G_1/N$. From Lemma \ref{lemreunit} it follows that $res_{G_1/N}(\varepsilon)=1$. Hence, for $\sigma,\tau\in G_1$, we have:
	\begin{displaymath}
		res_{G_1}(e)(\sigma,\tau)=e(\sigma,\tau)=\varepsilon(\sigma N,\tau N)=res_{G_1/N}(\varepsilon)(\sigma N,\tau N)=1,
	\end{displaymath}
since $\sigma N,\tau N\in G_1/N$. Then $res_{G_1}(e)=1$ and again from Lemma \ref{lemreunit} we get $G_1\leq H$ from where $N\leq H$.

Finally we will show that $G_1=H$ and then the equality of the statement of the theorem is obvious. If $\sigma\in H$, then $e(\sigma,\sigma^{-1})\neq 0$. Since $N\trianglelefteq G$ it holds that $(\sigma N)^{-1}=\sigma^{-1}N$ and so
	\begin{displaymath}
		\varepsilon(\sigma N,(\sigma N)^{-1})=\varepsilon(\sigma N,\sigma^{-1}N)=e(\sigma,\sigma^{-1})\neq 0.
	\end{displaymath}
But then, $\sigma N\in H_{G/N}(\varepsilon)=G_1/N$ and so $\sigma\in G_1$. Hence $G_1=H$.\quad $\square$\\

A special case of Theorem \ref{theinertquat} is when $N$ equals the inertial group of $e$. Then $e$ comes always from a lifting of an idempotent $2$-cocycle with trivial inertial group, as the next corollary shows.

\begin{corollary}{\label{corinflation}}
Let $e\in E^2(G,L;H)$. If $H$ is a normal subgroup of $G$, then there exists $\varepsilon\in E^2(G/H,L^H;1)$ such that $e=\widehat{\varepsilon}$.
\end{corollary}
\textbf{Proof}: From Theorem \ref{theinertquat}, setting $N=H$, it follows that there exists $\varepsilon\in E^2(G/H,L^H)$ such that $e=\widehat{\varepsilon}$. Furthermore $H_{G/H}(\varepsilon)=H_G(\widehat{\varepsilon})/H=H_G(e)/H=1$.\quad$\square$\\

The opposite is also true.

\begin{corollary}{\label{corinftrweak}}
Let $H\unlhd G$ and $\varepsilon\in E^2(G/H,L^H;1)$. Then there exists $e\in E^2(G,L;H)$ such that $e=\widehat{\varepsilon}$.
\end{corollary}
\textbf{Proof}: By inflating $\varepsilon$ from $G/H$ to $G$, we construct $e\in E^2(G,L)$ such that $e=\widehat{\epsilon}$. We only have to show that the inertial group of $e$ is $H$. Indeed
	\begin{eqnarray}
		H_G(e) &=& \{\sigma\in G:e(\sigma,\sigma^{-1})=1\}=\{\sigma\in G:\varepsilon(\sigma H,\sigma^{-1}H)=1\}\nonumber\\
			&=& \{\sigma\in G: \sigma H\in H_{G/H}(\varepsilon)\}=H.\nonumber\quad\square
	\end{eqnarray}

Let $H\unlhd G$ and $r\in Sl(G/H)$. In this case, lifting the idempotent $2$-cocycle $\varepsilon_r$ from $G/H$ to $G$ can also be accomplished using (\cite{LaTh}, Proposition 8.6 and Corollary 8.7). Specifically, if we set $\widehat{r}(\sigma)=r(\sigma H)$, for $\sigma\in G$, then $\widehat{r}\in Sl(G)$ with $M_{\widehat{r}}=H$. For the corresponding idempotent $2$-cocycles we have the relation $e_{\widehat{r}}=\widehat{\varepsilon}_r$.

\newpage

\begin{example}{\label{examliftdih}}
\end{example}
Let
	\begin{displaymath}
		G=D_8=<\sigma,a:a^4=\sigma^2=1,\sigma a\sigma^{-1}=a^{-1}>
	\end{displaymath}
be the dihedral group of order $8$ and identity element $1$. Let $H=\{1,a^2\}$. It is a normal subgroup of $G$ with $G/H=\{H,aH,\sigma H,a\sigma H\}$. In order to construct an element $r$ of $Sl(G/H)$ with trivial inertial group, we start with a generating set of $(G/H)^*=G/H\setminus \{H\}$ as in (\cite{LaTh},Section 10), say $A_1=\{aH,\sigma H\}$. Then $A_2=\{a\sigma H\}$ and $r$ takes the following values:
	\begin{displaymath}
		r(H)=0,\;r(aH)=1,\;r(\sigma H)=1,\;r(a\sigma H)=2.
	\end{displaymath}

\noindent The corresponding idempotent $2$-cocycle $\varepsilon_r:G/H\times G/H\longrightarrow L^H$ takes the following values:
	\begin{displaymath}
		\begin{array}{|r|cccc|}
		\hline
			(\sigma,\tau) & H & aH & \sigma H & a\sigma H\\
			\hline
			H    & 1 & 1 & 1 & 1\\
			aH   & 1 & 0 & 1 & 0\\
			\sigma H & 1 & 1 & 0 & 0\\
			a\sigma H & 1 & 0 & 0 & 0\\
			\hline
		\end{array}
	\end{displaymath}
As in (\cite{LaTh},Proposition 8.6) we construct $\widehat{r}\in Sl(G)$ with the following values:
 	\begin{displaymath}
		\begin{array}{ll}
			\widehat{r}(1)=\widehat{r}(a^2)=r(H)=0, & \widehat{r}(a)=\widehat{r}(a^3)=r(aH)=1,\\
			\widehat{r}(\sigma)=\widehat{r}(a^2\sigma)=r(\sigma H)=1, & \widehat{r}(a\sigma)=\widehat{r}(a^3\sigma)=r(a\sigma H)=2.
		\end{array}
	\end{displaymath}
>From the relation (\ref{eqidcocyle}) we construct the $2$-cocyle $e_{\widehat{r}}:G\times G\longrightarrow L$ with the following values:
	\begin{displaymath}
		\begin{array}{c|c|cc|cc|cc|cc|}
			 &        & H & & aH & & \sigma H & & a\sigma H & \\
		\hline
			 &(\sigma,\tau)   & 1 & a^2 & a & a^3 & \sigma & a^2\sigma & a\sigma & a^3\sigma\\
			\hline
			H & 1     & 1&1&1&1&1&1&1&1\\
			  &a^2    & 1&1&1&1&1&1&1&1\\
			\hline
			aH & a    & 1&1&0&0&1&1&0&0\\
			   & a^3  & 1&1&0&0&1&1&0&0\\
			\hline
			\sigma H & \sigma    & 1&1&1&1&0&0&0&0\\
			   & a^2\sigma & 1&1&1&1&0&0&0&0\\
			\hline
			a\sigma H& a\sigma   & 1&1&0&0&0&0&0&0\\
			   & a^3\sigma & 1&1&0&0&0&0&0&0\\
			\hline
		\end{array}
	\end{displaymath}

\noindent The corresponding graphs on left cosets are:
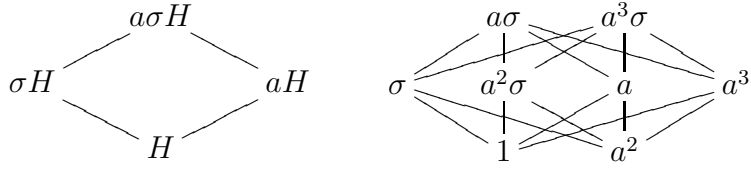
\[
\xymatrixrowsep{0.15in}
\xymatrixcolsep{0.3in}
\xymatrix{ 	& a\sigma H &\\
					\sigma H\ar@{-}[ur]& & aH\ar@{-}[ul]\\
					& H \ar@{-}[ul]\ar@{-}[ur] &}\qquad
			\xymatrix{ 	& a\sigma  &  a^3\sigma & \\
					\sigma \ar@{-}[ur]\ar@{-}[urr] &  a^2\sigma\ar@{-}[u]\ar@{-}[ur] & a \ar@{-}[ul]\ar@{-}[u] & a^3\ar@{-}[ull]\ar@{-}[ul]\\
					& 1 \ar@{-}[ul]\ar@{-}[u]\ar@{-}[ur]\ar@{-}[urr] & a^2\ar@{-}[ull]\ar@{-}[ul]\ar@{-}[u]\ar@{-}[ur] & } 		 
\]

\noindent Finally, the values of the function $\epsilon_{\widehat{r}}:G\times G\longrightarrow S^*$ from relation (\ref{eqepsilonr}) are the following:
	\begin{displaymath}
		\begin{array}{c|c|cc|cc|cc|cc|}
			 &        & H & & aH & & \sigma H & & a\sigma H & \\
		\hline
			 &(\sigma,\tau)   & 1 & a^2 & a & a^3 & \sigma & a^2\sigma & a\sigma & a^3\sigma\\
			\hline
			H & 1     & 1&1&1&1&1&1&1&1\\
			  &a^2    & 1&1&1&1&1&1&1&1\\
			\hline
			aH & a    & 1&1&\pi^2&\pi^2&1&1&\pi^2&\pi^2\\
			   & a^3  & 1&1&\pi^2&\pi^2&1&1&\pi^2&\pi^2\\
			\hline
			\sigma H & \sigma    & 1&1&1&1&\pi^2&\pi^2&\pi^2&\pi^2\\
			   & a^2\sigma & 1&1&1&1&\pi^2&\pi^2&\pi^2&\pi^2\\
			\hline
			a\sigma H& a\sigma   & 1&1&\pi^2&\pi^2&\pi^2&\pi^2&\pi^4&\pi^4\\
			   & a^3\sigma & 1&1&\pi^2&\pi^2&\pi^2&\pi^2&\pi^4&\pi^4\\
			\hline
		\end{array}
	\end{displaymath}

\noindent We remind that $\epsilon_{\widehat{r}}$ is not a $2$-cocycle. In order $\epsilon_{\widehat{r}}$ to become an element of $B^2(G,S^*)$, every element $\epsilon_{\widehat{r}}(\sigma,\tau)$, with $\sigma,\tau\in G$, must be multiplied with the invertible term $u(\sigma)^{r(\tau)}$.\\

\end{document}